# Exploring two concepts: branch decomposition and weak ultrafilter on connectivity system

Takaaki Fujita
Affiliation:Independent
Email:t171d603@gunma-u.ac.jp

**Abstract:** This paper explores two fundamental concepts: branch width and weak ultrafilter. Branch width is a significant graph width parameter that measures the degree of connectivity in a graph using a branch decomposition and a symmetric submodular function. Weak ultrafilter, introduced as a weakened definition of an ultrafilter, plays a vital role in interpreting defaults in logic. We introduce the concept of Weak Ultrafilter on the connectivity system *(X, f)* and demonstrate its duality with branch decomposition. This study enhances our understanding of these concepts in graph combinatorial and logical contexts.

## 1. Introduction

In this paper, we aim to explore two important concepts: branch width and weak ultrafilter. Let's begin with branch width. In graph theory, branch width is one of important graph width parameter. It involves a branch decomposition, where the leaves of the decomposition correspond to the edges of a graph. Each edge is associated with a value from a symmetric submodular function, which measures the degree of connectivity between the edges. Branch width offers the advantage of generalizing the width of symmetric submodular functions on graphs. Its study has been of great importance in graph theory due to its significant impact on the field. Various graph width parameters, such as path-width, tree-width, carving-width, cut-width, and linear-width, have been extensively studied in addition to branch-width. Due to their significance, numerous papers have been published on these topics (e.g., [6-25]).

Now let's turn to the concept of weak ultrafilter. In the realm of logic, weak filter (weak ultrafilter) is a concept introduced by K. Schlechta [1] in the 1990s as a weakened definition of an ultrafilter. It serves as a powerful tool for interpreting defaults through a generalized "most" quantifier in first-order logic. Weak ultrafilters can be considered the broadest class of structures that qualify as a "collection of large subsets" of a given index set. They are difficult to define in a weaker, yet still plausible, manner. The notion of weak ultrafilter naturally emerges and has found applications in epistemic logic and other areas of knowledge representation. Its study has been the subject of numerous research papers (ex. [1,2,3,4,5] ), highlighting its importance in the field.

In this short paper, we introduce the concept of a Weak Ultrafilter on the connectivity system *(X, f),* and we demonstrate its duality with branch decomposition, a graph width parameter. The connectivity system *(X, f)* refers to a pair consisting of an underlying set *X* (assumed to be the set representing the entire field of discussion) and a symmetric submodular function f. By examining the connectivity system, we can gain insights into the properties and relationships between weak ultrafilters and branch decomposition.

Overall, the exploration of branch width and weak ultrafilter in this paper contributes to a deeper understanding of these concepts and their relevance in various mathematical and logical contexts.

## 2. Definitions and Notations in this paper

In this section, we present the mathematical definitions and notations for each concept.

### 2.1 Symmetric Submodular Function and Connectivity System

The definition of a symmetric submodular function is given below. The symmetric submodular function is widely utilized and discussed in various scholarly publications (e.g., [41-44]).

**Definition 1 :** Let $X$ be a finite set. A function $f: X \to \mathbb{N}$ is called symmetric submodular if it satisfies the following conditions:
- $\forall A \subseteq X, f(A) = f(X \setminus A)$.
- $\forall A, B \subseteq X, f(A) + f(B) \geq f(A \cap B) + f(A \cup B)$.

In this short paper, a pair $(X, f)$ of a finite set $X$ and a symmetric submodular function $f$ is called a connectivity system.
It is known that a symmetric submodular function $f$ satisfies the following properties:
**Lemma 2[11] :** A symmetric submodular function $f$ satisfies:
1. $\forall A \subseteq X, f(A) \geq f(\emptyset) = f(X)$.
2. $\forall A, B \subseteq X, f(A) + f(B) \geq f(A \setminus B) + f(B \setminus A)$.

In this paper, we use the notation $f$ for a symmetric submodular function, a finite set $X$, and a natural number $k$. A set A is $k$-efficient if $f(A) \leq k$. Unless otherwise specified, in this paper, the underlying set $X$ is assumed to be a non-empty finite set.

## 2.2 Weak Ultrafilter on Boolean Algebra
The definition of a weak ultrafilter in a Boolean algebra $(X, \cup, \cap)$ is given below. Note that in the realm of logic, weak filter (weak ultrafilter) is a concept as a weakened definition of an ultrafilter.

**Definition 3[16]:** In a Boolean algebra $(X, \cup, \cap)$, a set family $F \subseteq 2^X$ satisfying the following conditions is called a weak filter on the carrier set X.
(F1) $A, B \in F \Rightarrow A \cap B \neq \emptyset$
(F2) $A \in F, A \subseteq B \subseteq X \Rightarrow B \in F$.
(F3) $\emptyset$ is not belong to $F$.

In a Boolean algebra $(X, \cup, \cap)$, A weak ultrafilter is called an weak ultrafilter and satisfies the following axiom (F4):
(F4) $\forall A \subseteq X$, either $A \in F$ or $X / A \in F$.

Ultrafilters, commonly defined on Boolean algebras as satisfying the axiom (F1): $A, B \in F \Rightarrow A \cap B \in F$, have broad and significant applications across various disciplines (ex.[27-37] ). They are extensively utilized in topology, algebra, logic, set theory, lattice theory, matroid theory, graph theory, combinatorics, measure theory, model theory, and functional analysis. Ultrafilters play a crucial role in comprehending mathematical structures, investigating infinite objects, and offering powerful tools for reasoning and analysis in diverse domains. Consequently, the study of Weak Ultrafilters, a concept closely related to Ultrafilters, holds utmost importance.

## 2.3 Weak Ultrafilter on connectivity system
The definition of Weak Ultrafilter is modified by incorporating additional conditions from Submodular functions. Specifically, Weak Ultrafilter on Boolean Algebra now includes the requirements related to Submodular functions.

**Definition 4:** In a connectivity system, the set family $W \subseteq 2^X$ is called a weak ultrafilter of order $k+1$ if the following axioms hold true:
(FB) For every $A \in W, f(A) \leq k$.
(FH) If $A, B \subseteq X, f(B) \leq k,$ A is a proper subset of $B$ and $A$ belongs to $W$, then $B$ belongs to $W$.

(WIS) If *A* belongs to *W*, B belongs to *W* and $f(A \cap B) \leq k$, then $A \cap B \neq \emptyset$.
(FW) $\emptyset$ does not belong to *W*.
(FE) If $\forall A \subseteq X$ and $f(A) \leq k$, then either $A \in W$ or $X / A \in W$.

And weak ultrafilter is non-principal if the weak ultrafilter satisfies following axiom:
(FP) $A \notin W$ for all $A \subseteq X$ with $|A| = 1$.

Note that Ultrafilters, commonly defined on connectivity system as satisfying the axiom (FS): If *A* belongs to *W*, *B* belongs to *W* and $f(A \cap B) \leq k$, then $A \cap B \notin W$.

### 2.4 Branch-decomposition of a connectivity system
The definition of branch-decomposition is shown below. Due to its significance, branch-decomposition has been the subject of various research studies [6,7,21,22,37-40].

**Definition 5:** Let *(X, f)* be a connectivity system. The pair *(T, μ)* is a branch decomposition tree of *(X, f)* if *T* is a ternary tree such that $|L(T)| = |X|$ and μ is a bijection from *L(T)* to *X*, where *L(T)* denotes the leaves in *T*. For each $e \in E(T)$, we define *bw(T, μ, e)* as $f(\cup_{v \in L(T1)} \mu(v))$, where $T_1$ is a tree obtained by removing *e* from *T* (taking into account the symmetry property of f). The width of *(T, μ)* is defined as the maximum value among *bw(T, μ, e)* for all $e \in E(T)$. The branch-width of *X*, denoted by *bw(X)*, is defined as the minimum width among all possible branch decomposition trees of *X*.

### 3. Characterization of Weak Ultrafilter
In this section, we aim to prove the dual theorems of branch width and Weak Ultrafilter as presented below.

**Theorem 6:** If *W* is a Weak Ultrafilter of order *k+1*, then branch width of the connectivity system *(X, f)* is at most *k*.

**Proof:** Let *(X, f)* be a connectivity system and *W* be a Weak Ultrafilter of order *k+1*. We think about $I = \{A \mid X \setminus A \in W\}$. Let *(T, μ)* be a branch decomposition tree of *(X, f)* with width *bw(T, μ)*.
For any edge *e* in *E(T)*, we have the branch-width as:
*bw(T, μ, e)* = $f(\cup_{v \in L(T1)} \mu(v))$ where $T_1$ is the tree obtained by removing *e* from *T*.

Now, as per Definition 4 of the Weak Ultrafilter, we have: For every $A \in I$, $f(A) \leq k$.
Therefore, for any set A in *I*, the function value *f(A)* is at most *k*.
Now consider the set $S = \cup_{v \in L(T1)} \mu(v)$ for some edge *e* in *E(T)*. As per the definition of Weak Ultrafilter (axiom (FE)), either *S* is in *I* or its complement is in *I*. If *S* is in *I*, *f(S)* is at most *k*, by axiom (IB). If the complement of *S* is in *I*, then $f(S) = f(X \setminus S)$, which is also at most *k* due to the symmetric property of *f*.
Therefore, *bw(T, μ, e)* is at most *k* for all *e* in *E(T)*, and hence the width *bw(T, μ)* of the branch decomposition is at most *k*. Since this holds for any branch decomposition *(T, μ)*, it follows that the branch-width *bw(X)* of *X* is also at most *k*. This proof is completed. ∎

**Theorem 7:** If branch-width of the connectivity system *(X, f)* is *k+1*, then no Weak Ultrafilter of order *k+1* exists.
**Proof:** Assume for contradiction that there exists a Weak Ultrafilter *W* of order *k+1* for a connectivity system *(X, f)* whose branch-width is *k+1*. Note that We think about $I = \{A \mid X \setminus A \in W\}$.
Let *(T, μ)* be a branch decomposition of the connectivity system *(X, f)* that achieves the branch-width *k+1*, i.e., the width of *(T, μ)* is *k+1*. This means there exists an edge *e* in *E(T)* such that:

$bw(T, \mu, e) = f(\bigcup_{v \in L(T_1)} \mu(v)) = k+1$

where $T_1$ is the tree obtained by removing $e$ from $T$.

Now, consider the set $S = \bigcup_{v \in L(T_1)} \mu(v)$. By axiom (FE) of the Weak Ultrafilter, either $S$ is in $I$ or its complement is in $I$. If $S$ is in $I$, then by axiom (GB), $f(S)$ should be at most $k$, which is a contradiction. Similarly, if the complement of $S$ is in $I$, then $f(S) = f(X \setminus S)$ should be at most $k$ due to the symmetric property of $f$, which is again a contradiction.

Hence, no Weak Ultrafilter of order $k+1$ can exist if the branch-width of the connectivity system $(X, f)$ is $k+1$. This proof is completed. ∎

## Acknowledgments

I humbly express my sincere gratitude to all those who have extended their invaluable support, enabling me to successfully accomplish this paper.

## Conflict of Interest Statement

The author declares no conflicts of interest.